\documentclass[12pt]{article}
\usepackage{amssymb, amsfonts, amsmath, amsthm}

\newtheorem{theorem}{Theorem}

\newtheorem{corollary}[theorem]{Corollary}

\newcommand{\calg}{{\cal G}}
\newcommand{\dd}{\displaystyle }

\def\3{\subset }
\def\4{\subseteq }

\def\0{\leqno}

\def\barr{\begin{array}}
\def\earr{\end{array}}
\def\dd{\displaystyle}
\def\Z{{\rlap{$\kern2pt{\rm Z}$}{\rm Z}\,}}

\def\calg{{\cal G}}

\title{\bf The normal subgroup structure\\ of ZM-groups}
\author{Marius T\u arn\u auceanu}
\date{February 17, 2015}

\begin{document}

\maketitle

\begin{abstract}
The main goal of this note is to determine and to count the normal
subgroups of a ZM-group. We also indicate some necessary and
sufficient conditions such that the normal subgroups of a ZM-group
form a chain.
\end{abstract}

\noindent{\bf MSC (2010):} Primary 20D30; Secondary 20D60, 20E99.

\noindent{\bf Key words:} ZM-groups, normal subgroups, chains.

\section{Introduction}

The starting point for our discussion in given by the paper
\cite{2}, where the class $\calg$ of finite groups that can be
seen as cyclic extensions of cyclic groups has been considered.
The main theorem of \cite{2} furnishes an explicit formula for the
number of subgroups of a group contained in $\calg$. In
particular, this number is computed for several remarkable
subclasses of $\calg$: abelian groups of the form $\mathbb{Z}_m
\times \mathbb{Z}_n$, dihedral groups $D_{2m}$, and Zassenhaus
metacyclic groups (ZM-groups, in short).

In group theory the study of the normal subgroups of (finite)
groups plays a very important role. So, the following question
concerning the class $\calg$ is natural:

{\it \hspace{10mm}Which is the number of {\it normal} subgroups of
a group in $\calg$?}

\noindent The purpose of the current note is to answer partially
this question, by finding this number for the above three
subclasses of $\calg$. Since all subgroups of an abelian group are
normal, for the first subclass the answer is given by \cite{2}.
The number of normal subgroups of the dihedral group $D_{2m}$ is
also well-known, namely $\tau(m)+1$ if $m$ is odd, and $\tau(m)+3$
if $m$ is even (as usually, $\tau(m)$ denotes the number of
distinct divisors of $m\in\mathbb{N}^*$). Therefore we will focus
only on describing and counting the normal subgroups of ZM-groups.

Most of our notation is standard and will not be repeated here.
Basic definitions and results on group theory can be found in
\cite{5,6,8}. For subgroup lattice theory we refer the reader to
\cite{7,9}.
\bigskip

First of all, we recall that a ZM-group is a finite group with all
Sylow subgroups cyclic. By \cite{5}, such a group is of type
\begin{equation}
{\rm ZM}(m,n,r)=\langle a, b \mid a^m = b^n = 1,
\hspace{1mm}b^{-1} a b = a^r\rangle, \nonumber
\end{equation}
where the triple $(m,n,r)$ satisfies the conditions
\begin{equation}
{\rm gcd}(m,n) = {\rm gcd}(m, r-1) = 1 \quad \text{and} \quad r^n
\equiv 1 \hspace{1mm}({\rm mod}\hspace{1mm}m). \nonumber
\end{equation}
It is clear that $|{\rm ZM}(m,n,r)|=mn$, ${\rm
ZM}(m,n,r)'\hspace{1mm}=\hspace{1mm}\langle a \rangle$
(consequently, we have $|{\rm ZM}(m,n,r)'|=m$) and ${\rm
ZM}(m,n,r)/{\rm ZM}(m,n,r)'$ is cyclic of order $n$. One of the
most important (lattice theoretical) property of the ZM-groups is
that these groups are exactly the finite groups whose poset of
conjugacy classes of subgroups forms a distributive lattice (see
Theorem A of \cite{1}). We infer that they are DLN-groups, that is
groups with distributive lattice of normal subgroups.
\bigskip

The subgroups of ${\rm ZM}(m,n,r)$ have been completely described
in \cite{2}. Set
\begin{equation}
L=\left\{(m_1,n_1,s)\in\mathbb{N}^3 \hspace{1mm}\mid\hspace{1mm}
m_1|m,\hspace{1mm} n_1|n,\hspace{1mm} s<m_1,\hspace{1mm}
m_1|s\frac{r^n-1}{r^{n_1}-1}\right\}. \nonumber
\end{equation}
Then there is a bijection between $L$ and the subgroup lattice
$L({\rm ZM}(m,n,r))$ of ${\rm ZM}(m,n,r)$, namely the function
that maps a triple $(m_1,n_1,s)\in L$ into the subgroup
$H_{(m_1,n_1,s)}$ defined by
\begin{equation}
H_{(m_1,n_1,s)}=\bigcup_{k=1}^{\frac{n}{n_1}}\alpha(n_1,
s)^k\langle a^{m_1}\rangle=\langle a^{m_1},\alpha(n_1, s)\rangle,
\nonumber
\end{equation}
where $\alpha(x, y)=b^xa^y$, for all $0\leq x<n$ and $0\leq y<m$.
Remark also that $|H_{(m_1,n_1,s)}|=\frac{mn}{m_1n_1}$, for any
$s$ satisfying $(m_1,n_1,s)\in L$.
\bigskip

By using this result, we are able to describe the normal subgroup
structure of ${\rm ZM}(m,n,r)$.

\begin{theorem}
The normal subgroup lattice $N({\rm ZM}(m,n,r))$ of ${\rm
ZM}(m,n,r)$ consists of all subgroups
\begin{equation}
H_{(m_1,n_1,s)}\in L({\rm ZM}(m,n,r))\hspace{2mm} {\rm with}
\hspace{2mm}(m_1,n_1,s)\in L', \nonumber
\end{equation}
where
\begin{equation}
L'=\left\{(m_1,n_1,s)\in\mathbb{N}^3 \hspace{1mm}\mid\hspace{1mm}
m_1|{\rm gcd}(m,r^{n_1}-1),\hspace{1mm} n_1| n,\hspace{1mm}
s=0\right\}\subseteq L. \nonumber
\end{equation}
\end{theorem}

We infer that, for every $m_1|m$ and $n_1|n$, ${\rm ZM}(m,n,r)$
possesses at most one normal subgroup of order
$\frac{mn}{m_1n_1}$. In this way, all normal subgroups of ${\rm
ZM}(m,n,r)$ are characteristic. In particular, the above theorem
allows us to count them.

\begin{corollary}
The following equality holds
\begin{equation}
|N({\rm ZM}(m,n,r))|=\dd\sum_{n_1\mid n}\tau({\rm
gcd}(m,r^{n_1}-1)). \tag{1}
\end{equation}
\end{corollary}

In the following we will denote by $d$ the multiplicative order of
$r$ modulo $m$, that is
$$d={\rm min}\hspace{1mm}\{k\in\mathbb{N}^* \mid r^k\equiv1 \hspace{1mm}({\rm mod} \hspace{1mm}m)\}.$$
Clearly, the sum in the right side of (1) depends on $d$. For $m$
or $n$ primes, this sum can be easily computed.

\begin{corollary}
If $m$ is a prime, then
\begin{equation}
|N({\rm ZM}(m,n,r))|=\tau(n)+\tau(\frac{n}{d}), \tag{2}
\end{equation}
while if $n$ is a prime, then
\begin{equation}
|N({\rm ZM}(m,n,r))|=\tau(m)+1. \tag{3}
\end{equation}
\end{corollary}

Mention that the number of normal subgroups of the dihedral group
$D_{2m}$ with $m$ odd can be obtained from (3), by taking $n=2$.
\bigskip

Next we will focus on finding the triples $(m,n,r)$ for which
$N({\rm ZM}(m,n,r))$ becomes a chain.

\begin{theorem}
The normal subgroup lattice $N({\rm ZM}(m,n,r))$ of ${\rm
ZM}(m,n,r)$ is a chain if and only if either $m=1$ and $n$ is a
prime power, or both $m$ and $n$ are prime powers and ${\rm
gcd}(m,r^k-1)=1$ for all $1\leq k<n$.
\end{theorem}
\bigskip

Remark that Theorem 4 gives a method to construct finite (both
abelian and nonabelian) groups whose lattices of normal subgroups
are chains of prescribed lengths.
\bigskip

Finally, we indicate an open problem with respect to the above
results.

\bigskip\noindent{\bf Open problem.} Describe and count the normal
subgroups of an {\it arbitrary} finite group contained in $\calg$.
Also, extend these problems to {\it arbitrary} finite metacyclic
groups, whose structure is well-known (see, for example,
\cite{4}).

\section{Proofs of the main results}

\bigskip\noindent{\bf Proof of Theorem 1.}
First of all, we observe that under the notation in Section 1 we
have
\begin{equation}
\alpha(x_1, y_1)\alpha(x_2, y_2)=\alpha(x_1+x_2, r^{x_2}y_1+y_2).
\nonumber
\end{equation}
This implies that
\begin{equation}
\alpha(x, y)^k=b^{kx}a^{y\frac{r^{kx}-1}{r^x-1}}, \text{ for all }
k\in\mathbb{Z}, \hspace{2mm}{\rm and}\hspace{2mm} \alpha(x,
y)^{-1}=\alpha(-x, -r^{-x}y). \nonumber
\end{equation}
Since
\begin{equation}
\alpha(x, y)^{-1}\alpha(n_1, s)\alpha(x, y)=\alpha(n_1, t_{x,y}),
\hspace{2mm} {\rm where} \hspace{2mm}t_{x,y}=-r^{n_1}y+r^xs+y,
\nonumber
\end{equation}
one obtains
\bigskip
$$\hspace{-40mm}H_{(m_1,n_1,s)}^{\alpha(x, y)}\hspace{-1mm}=\alpha(x, y)^{-1}H_{(m_1,n_1,s)}\alpha(x, y)\hspace{-1mm}=$$
$$\hspace{-3mm}=\bigcup_{k=1}^{\frac{n}{n_1}}\alpha(x, y)^{-1}\alpha(n_1,
s)^k\alpha(x,y)^{-1}\langle a^{m_1}\rangle=$$
$$\hspace{-2mm}=\bigcup_{k=1}^{\frac{n}{n_1}}\left(\alpha(x, y)^{-1}\alpha(n_1,
s)\alpha(x, y)\right)^k\langle a^{m_1}\rangle=$$
$$\hspace{-9mm}=\bigcup_{k=1}^{\frac{n}{n_1}}\alpha(n_1, t_{x,y})^k\langle a^{m_1}\rangle=H_{(m_1,n_1,t_{x,y})}$$with
the convention that $t_{x,y}$ is possibly replaced by $t_{x,y}
\hspace{1mm}{\rm mod} \hspace{1mm}m_1$. Then $H_{(m_1,n_1,s)}$ is
normal in ${\rm ZM}(m,n,r)$ if and only if we have $t_{x,y}\equiv
s \hspace{1mm}({\rm mod}\hspace{1mm} m_1)$, or equivalently
\begin{equation}
m_1|s(r^x-1)-y(r^{n_1}-1), \tag{4}
\end{equation}
for all $0\leq x<n$ and $0\leq y<m$. Take $x=0$ in (4). It follows
that $m_1|y(r^{n_1}-1)$, for all $0\leq y<m$, and so
$m_1|r^{n_1}-1.$ We get $m_1|s(r^x-1)$, for all $0\leq x<n$. By
putting $x=1$ and using the equality gcd($m,r-1$)=1, it results
$m_1|s$. But $s<m_1$, therefore $s=0$. Hence we have proved that
the subgroup $H_{(m_1,n_1,s)}$ is normal if and only if $m_1|{\rm
gcd}(m,r^{n_1}-1)$ and $s=0$, as desired.
\hfill\rule{1,5mm}{1,5mm}

\bigskip\noindent{\bf Proof of Theorem 4.}
Suppose first that $N({\rm ZM}(m,n,r))$ is a chain. Then ${\rm
ZM}(m,n,r)$ is a monolithic group, that is it possesses a unique
minimal normal subgroup. By Theorem 5.9 of \cite{3} it follows
that either $m=1$ and $n$ is a prime power, or $m$ is a prime
power and $r^k \not\equiv 1 \hspace{1mm}({\rm mod}\hspace{1mm} m)$
for all $1 \leq k < n$. On the other hand, we observe that $N({\rm
ZM}(m,n,r))$ contains the sublattice
\begin{equation}
L_1=\left\{H_{(1,n_1,0)} \hspace{1mm}\mid \hspace{1mm}n_1|n
\right\}, \nonumber
\end{equation}
which is isomorphic to the lattice of all divisors of $n$. Thus
$n$ is a prime power, too. In order to prove the last assertion,
let us assume that
\begin{equation}
{\rm gcd}(m,r^k-1)=m_1\neq 1 \nonumber
\end{equation}
for some $1\leq k<n$ and consider $k$ to be minimal with this
property. It follows that $k|n$. Then the subgroup $H_{(m_1,k,0)}$
belongs to $N({\rm ZM}(m,n,r))$, but it is not comparable to
$H_{(1,n,0)}={\rm ZM}(m,n,r)'$, a contradiction.

Conversely, if the triple $(m,n,r)$ satisfies one of the
conditions in Theorem 4, then $N({\rm ZM}(m,n,r))$ is either a
chain of length $v$ for $m=1$ and $n=q^v$ ($q$ prime), namely
\begin{equation}
H_{(1,q^v,0)}\subset H_{(1,q^{v-1},0)}\subset \cdots \subset
H_{(1,1,0)}, \nonumber
\end{equation}
or a chain of length $u+v$, for $m=p^u$ and $n=q^v$ ($p,q$
primes), namely
\begin{equation}
H_{(p^u,q^v,0)}\subset H_{(p^{u-1},q^v,0)}\subset \cdots \subset
H_{(1,q^v,0)}\subset H_{(1,q^{v-1},0)}\subset \cdots \subset
H_{(1,1,0)}. \nonumber
\end{equation}
This completes the proof.
\hfill\rule{1,5mm}{1,5mm}

\bigskip\noindent{\bf Acknowledgements.} The author is grateful to the reviewer for
its remarks which improve the previous version of the paper.

\vspace*{5ex}\small

\hfill
\begin{minipage}[t]{5cm}
Marius T\u arn\u auceanu \\
Faculty of  Mathematics \\
``Al.I. Cuza'' University \\
Ia\c si, Romania \\
e-mail: {\tt tarnauc@uaic.ro}
\end{minipage}

\end{document}